\documentclass[10 pt,leqno]{amsart}

\usepackage {amsfonts}
\usepackage{amsthm}
\usepackage{amssymb}
\usepackage{latexsym}
\usepackage{amsmath}
\usepackage{mathrsfs}

\newenvironment{rlist}
{

\begin{enumerate}}
{\end{enumerate}}

\theoremstyle{plain}
\newtheorem{tw}{Theorem} [section]

\newtheorem {prop}[tw] {Proposition}

\theoremstyle{definition}
\newtheorem {deft}[tw] {Definition}

\newcommand{\bn}{\Bbb N}
\newcommand{\br}{\Bbb R}

\newcommand{\pred}{\mlg_*}
\newcommand{\starpred}{\mlg_*^{\#}}
\newcommand{\modpred}{\mlg_*^{\tau}}

\newcommand{\dens}{D}

\newcommand{\pprdens}{\dens^{\frac{1}{p'}}}
\newcommand{\pdens}{\dens^{\frac{1}{2p}}}
\newcommand{\prdens}{\dens^{\frac{1}{p}}}
\newcommand{\fourdens}{\dens^{\frac{1}{4}}}
\newcommand{\twodens}{\dens^{\frac{1}{2}}}

\newcommand{\alg} {\mathsf{A}}
\newcommand{\Com} { \Delta}
\newcommand{\mlg} {\mathsf{M}}

\newcommand {\Lin} {{\textrm{Lin}}}

\newcommand{\antip}{\mathcal{S}}
\newcommand{\Tomita}{\mathcal{T}}
\newcommand{\Alg}{\mathcal{A}}

\newcommand{\Hil}{\mathsf{H}}

\newcommand{\lla}{\left\langle}
\newcommand{\rra}{\right\rangle}

\newcommand{\ot}{\otimes}
\newcommand{\ol}{\overline}

\newcommand{\wot}{\ol{\otimes}}
\newcommand{\ida}{\textrm{id}_{\alg}}
\newcommand{\idm}{\textrm{id}_{\mlg}}

\numberwithin{equation}{section}


\begin{document}

\author{Uwe Franz}
\address{D\'epartement de math\'ematiques de Besan\c{c}on,
Universit\'e de Franche-Comt\'e, 16, route de Gray, F-25 030
Besan\c{c}on cedex, France}
\curraddr{Graduate School of Information
Sciences, Tohoku University, Sendai 980-8579, Japan}
\urladdr{http://www-math.univ-fcomte.fr/pp\underline{ }Annu/UFRANZ/}
\email{uwe.franz at univ-fcomte.fr}
\thanks{U.F.\ was supported by a Marie Curie Outgoing International
Fellowship of the EU (Contract Q-MALL MOIF-CT-2006-022137) and a Polonium cooperation.}

\author{Adam Skalski}

\footnote{\emph{Permanent address of the first named author:} Department of Mathematics, University of \L\'{o}d\'{z}, ul. Banacha
22, 90-238 \L\'{o}d\'{z}, Poland.}

\address{Department of Mathematics and Statistics,  Lancaster University,
Lancaster, LA1 4YF} \email{a.skalski at lancaster.ac.uk}

\title{\bf On ergodic properties of convolution operators associated with compact quantum groups}

\begin{abstract}
Recent results of M.\,Junge and Q.\,Xu on the ergodic properties of the averages of kernels in noncommutative
$L^p$-spaces are applied to the analysis of the almost uniform convergence of operators induced by the
convolutions on compact quantum groups.
\end{abstract}
\maketitle

\vspace*{0.5cm} The classical ergodic theory was initially concerned with investigating the limits of iterations (or
iterated averages) of certain transformations of a measure space. Corresponding limit theorems were very
quickly seen to have natural generalisations in terms of the evolutions induced by operators acting on the
associated $L^p$-spaces (for an excellent treatment we refer to \cite{Krengel}).
The noncommutative counterpart of this theory is concerned with
investigation of limit properties for the iterations of operators acting on von Neumann algebras (viewed as
generalisations of classical $L^{\infty}$-spaces) or, more generally, on noncommutative $L^p$-spaces associated with
a von Neumann algebra equipped with a faithful normal state. It turned out that, after introducing appropriate
counterparts of the classical notion of almost everywhere convergence, one may consider in this generalised
context not only mean ergodic theorems, but also `pointwise' ones. This has been investigated intensively in the
70s and 80s by C.E.\,Lance, F.\,Yeadon, R.\,Jajte and others. Several results were obtained
for both the evolutions on von Neumann algebras and on $L^p$-spaces associated with a faithful normal trace.
 Recently M.\,Junge and Q.\,Xu in a beautiful paper
\cite{full} (whose main results were earlier announced in \cite{Junge}) proved new noncommutative maximal
inequalities and thus extended many ergodic theorems to the context of  Haagerup $L^p$-spaces, which
naturally arise when the considered state is non-tracial.

In this paper we apply the results of \cite{full} to obtain ergodic theorems for the evolutions induced by the
convolution operators on compact quantum groups (\cite{wor1}).
Although it is generally natural to view compact quantum groups as $C^*$-algebras, due to the nature of the
problems considered we prefer the von Neumann algebraic framework. It arises naturally as every
compact quantum group is equipped with a Haar state and one can pass to the corresponding  GNS representation. The
importance of this approach, where the Haar functional is a central notion from which in a sense the whole theory
is developed, is fully revealed in the context of locally compact quantum groups (\cite{kuv}). Here it provides
us both with a von Neumann algebra and with a canonical reference state on it.

The plan of the paper is as follows: after establishing notation and quoting preliminary results in the first
section, in Section 2 we introduce the convolution operators and obtain the ergodic theorems for their actions on
a compact quantum group $\mlg$. Section 3 contains a   discussion of the extensions to the case of Haagerup
$L^p$-spaces associated with the Haar state on $\mlg$ and in Section 4 we signal possible directions of
further investigations.

\section{Notations and preliminary results}

The symbol $\ot$ will denote the spatial tensor product of $C^*$-algebras, $\wot$ the ultraweak tensor product of
von Neumann algebras (and relevant extension of the algebraic tensor product of normal maps); $\odot$ will be
reserved for the purely algebraic tensor product.

\subsection*{Compact quantum groups}

The notion of compact quantum groups has been introduced in \cite{wor1}. Here we adopt the definition from
\cite{wor2}:

\begin{deft}
A \emph{compact quantum group} is a pair $(\alg, \Com)$, where $\alg$ is a unital $C^*$-algebra,
 $\Com:\alg \to \alg \ot \alg$ is a unital, *-homomorphic map which is
coassociative:
\[ (\Com \ot \ida) \Com = (\ida \ot \Com) \Com\]
 and $\alg$ satisfies the quantum cancellation properties:
\[ \overline{\Lin}((1\ot \alg)\Com(\alg) ) = \overline{\Lin}((\alg \ot 1)\Com(\alg) )
= \alg \ot \alg. \]
\end{deft}

One of the most important features of compact quantum groups is the existence of the dense $^*$-subalgebra $\Alg$
(the algebra of matrix coefficients of irreducible unitary representations of $\alg$), which is in fact a Hopf
$^*$-algebra - so for example $\Com: \Alg \to \Alg \odot \Alg$. As explained in the introduction, for us it is
more convenient to work in the von Neumann algebraic context.

\begin{deft}
A \emph{von Neumann algebraic (vNa) compact quantum group} is a pair $(\mlg, \Com)$, where $\mlg$ is a von Neumann
algebra,
 $\Com:\mlg \to \mlg \wot \mlg$ is a normal unital, *-homomorphic map which is
coassociative:
\[ (\Com \wot \ida) \Com = (\ida \wot \Com) \Com\]
and there exists a faithful normal state $h \in \mlg_*$ (called a Haar state) such that for all $x \in M$
\[ (h \wot \idm)\circ  \Com (x) = ( \idm \wot h)\circ  \Com (x) = h(x) 1.\]
\end{deft}

The next lemma and the comments below it should help to understand the connection between these two types of
objects.

\begin{prop}(\cite{wor2})
Let $\alg$ be a compact quantum group. There exists a unique state $h \in \alg^*$ (called the \emph{Haar state} of
$\alg$) such that for all $a \in \alg$
\[ (h \ot \ida)\circ  \Com (a) = ( \ida \ot h)\circ  \Com (a) = h(a) 1.\]
\end{prop}

A compact quantum group is said to be \emph{in reduced form} if the Haar state $h$ is faithful. If it is not the
case we can always quotient out the null space of $h$ ($\{a\in \alg: h(a^*a) =0\}$). This procedure in particular
does not influence the underlying Hopf $^*$-algebra $\Alg$; in fact the reduced object may be viewed as the
natural completion of $\Alg$ in the GNS representation with respect to $h$ (as opposed for example to the
universal completion of $\Alg$, for details see \cite{BMT}). We will therefore always assume that our compact
quantum groups are in reduced forms.

Let $\alg$ be a  compact quantum group and let $(\pi_h, \Hil)$ be the (faithful) GNS representation with respect
to the Haar state of $\alg$. Define $\mlg = \pi_h(\alg)''$. Then $\mlg$ is a von Neumann algebra, the coproduct
has a normal extension to $\mlg$ (denoted further by the same symbol) with values in $\mlg \wot \mlg$ and by the construction the Haar state retains its invariance properties in this new framework - we obtain the vNa compact quantum group. Conversely, given a vNa compact quantum group there is a way of associating to it a $C^*$-algebraic
objext, which is a compact quantum group (see [KuV$_{1-2}$] for the details of this construction and the statements which follow). As applying these
constructions twice yields the same (i.e.\ isomorphic) object as the original one, we can without the loss of
generality assume that whenever a vNa compact quantum group $(\mlg,\Com)$ is considered, it is in its standard
form given by a GNS representation with respect to the Haar state and that it has a w$^*$-dense unital
$C^*$-subalgebra $\alg$ such that $(\alg, \Com|_{\alg})$ is a compact quantum group.

Whenever $(\mlg, \Com)$ is a vNa compact quantum group, there exists a $^*$-antiauto-morphism of $\mlg$ (called
the \emph{unitary antipode} and denoted by $R$) and a $\sigma$-strongly$^*$ continuous one parameter group $\tau$
of $^*$-automorphisms of $\mlg$ (called a \emph{scaling group} of $(\mlg, \Com)$) such that the set $\Lin\{(\idm
\wot h)\left(\Com(x)(1 \ot y)\right): x,y \in \mlg\}$ is contained in the domain of a (densely defined) operator
$\antip = R \tau_{-\frac{1}{2}}$, called the \emph{antipode}. In fact the above set is a $\sigma$-strong$^*$ core
for $\antip$ and
\[ \antip \left((\idm \wot h)\left(\Com(x)(1 \ot y)\right) \right) =
 (\idm \wot h)\left((1 \ot x)\Com(y)\right)\;\;\;  x,y \in \mlg.\]
The unimodularity of compact quantum groups is expressed by the condition $h = h \circ R$ - in general the unitary
antipode exchanges the left invariant and the right invariant weights. Therefore we also have (by the strong left
invariance of the antipode)
\[ \antip \left((h \wot \idm )\left(\Com(x)(1 \ot y)\right) \right) =
 (h \wot\idm )\left((1 \ot x)\Com(y)\right)\;\;\;  x,y \in \mlg.\]
Additionally denote by $\Tomita$ the algebra of all analytic elements with respect to the modular group
(\cite{opalg}).

The coassociativity of $\Com$ implies that the predual of $\mlg$ equipped with the convolution product
\[ \phi \star \psi = (\phi \wot \psi) \Com,\; \; \; \;\phi, \psi \in \mlg_*\]
is a Banach algebra. It contains an important dense subalgebra that may be equipped with the involution relevant
for considering noncommutative counterparts of symmetric measures. Define, following \cite{kuv},
\[\starpred = \{ \omega \in \pred: \exists_{\theta \in \pred} \; \theta(x) = \ol{\omega} (\antip(x))
\textrm{ for all } x \in D(\antip)\}.\] The involution $^*$ in $\starpred$ is introduced with the help of the
obvious formula: $ \omega^* \supset \ol{\omega} \circ \antip$.

The modular group of the Haar state will be denoted simply by $\sigma$. Let us gather here a few useful
commutation relations:
\begin{equation} \label{copmod} (\tau_t \ot \sigma_t ) \Com =
        (\sigma_t \ot \tau_{-t} ) \Com   =  \Com \circ \sigma_t, \end{equation}
\begin{equation} \label{coptau} (\tau_t \ot \tau_t ) \Com =
        \Com \circ \tau_t, \end{equation}
\begin{equation} \label{Rscal} R \circ \tau_t =  \tau_{t} \circ R. \end{equation}

\subsection*{Notions of `pointwise' convergence in the von Neumann algebraic context}

Let $\mlg$ be a von Neumann algebra with a faithful normal state $\phi \in \mlg_*$, called the \emph{reference
state}.

\begin{deft}         \label{auconv}
A sequence $(x_n)_{n=1}^{\infty}$ of operators of $\mlg$ is almost uniformly (a.u.) convergent to $x \in \mlg$ if
for each $\epsilon >0 $ there exists $e \in P_M$  such that $ \phi (e^{\perp}) < \epsilon$ and
\[ \left\| (x_n -x) e \right\|_{\infty} \stackrel{n \longrightarrow  \infty}
  {\longrightarrow} 0.\]
A sequence $(x_n)_{n=1}^{\infty}$ of operators in $\mlg$ is bilaterally almost uniformly (b.a.u.) convergent to $x
\in \mlg$ if for each $\epsilon >0 $ there exists $e \in P_M$ such that $ \phi (e^{\perp}) < \epsilon$ and
\[ \left\| e(x_n -x) e \right\|_{\infty} \stackrel{n \longrightarrow  \infty}
  {\longrightarrow} 0.\]
\end{deft}

\begin{deft}
A linear map $T : \mlg \to \mlg$ is called a kernel (or a positive $L^1 - L^{\infty}$ contraction) if it is a
positive contraction:
\[ \forall_ {x\in \mlg} { \; \; 0 \leq x \leq I \Longrightarrow 0 \leq T (x) \leq I}\]
and has the property
\[ \forall_ {x\in \mlg, \, x\geq 0} {\; \; \phi (T(x)) \leq \phi (x) }.\]
\end{deft}

It is well known that for each kernel $T$ and $x \in \mlg$ the sequence $ (M_n(T) (x))_{n=1}^{\infty}$, where
\begin{equation} \label{aver} M_n (T) (x) = \frac{1}{n}\sum_{k=1}^{n} T^k (x), \end{equation}
is $w^*$-convergent to $F(x)$, where $F:\mlg \to \mlg$ denotes the w$^*$-continuous projection on the space of
fixed points of $T$.

The following individual ergodic theorem is due to B.\,K\"ummerer :

\begin{tw}[\cite{Kumm}] \label{ind}
 If $T:\mlg \to \mlg$ is a kernel, then for each $x\in \mlg$ the sequence $ (M_n(T) (x))_{n=1}^{\infty}$ converges to $F(x)$ almost uniformly.
\end{tw}

\section{Convolution operators and ergodic theorems on the level of a von Neumann algebra}

Let $(\mlg, \Com)$ be a vNa compact quantum group with the Haar state $h\in \mlg_*$. For any  $\phi \in \mlg_*$ by
the convolution operator associated with $\phi$ we shall understand the map $T_{\phi}: \mlg \to \mlg$ defined by
\begin{equation}  T_{\phi}  = (\idm \wot \phi)  \Com. \end{equation}
There is also an obvious left version, given by
\begin{equation}  L_{\phi}  = (\phi \wot \idm)  \Com. \end{equation}

The basic properties of the convolution operators are summarised below:

\begin{prop} \label{prr}
Let $\phi, \phi_i  \in \mlg_* (i \in \mathcal{I})$. Then the following hold:
\begin{rlist}
\item if $\phi \in \mlg^+_*$ then $T_{\phi}$ is completely positive; if $\phi(1)=1$ then $T_{\phi}$ is unital;
\item $T_{\phi}$ is normal and decomposable (the latter means it can be represented as a linear combination of completely
positive maps);
\item the map $\phi \longrightarrow T_{\phi}$ is a contractive homomorphism between Banach algebras $\mlg_*$ and $B(\mlg)$;
\item $ h \circ T_{\phi} = \phi(1) h$;
\item if $\phi_i \stackrel{i \in \mathcal{I}}{\longrightarrow} \phi $ in norm then $T_{\phi_i} \stackrel{i \in \mathcal{I}}{\longrightarrow} T_{\phi} $
in norm;
\item if $\phi_i \stackrel{i \in \mathcal{I}}{\longrightarrow} \phi $ weakly then for each $x \in \mlg$ $T_{\phi_i} (x)\stackrel{i \in \mathcal{I}}{\longrightarrow} T_{\phi} (x) $
in w$^*$-topology.
\end{rlist}
\end{prop}

\begin{proof}
Property (i) is obvious (as positive functionals are automatically CP), (ii) follows from (i) and the existence of
Jordan decomposition of normal functionals. Property (iii) is a consequence of coassociativity, contractivity of
$\Com$  and the fact that for each linear functional the completely bounded norm is equal to the standard norm.
(iv) follows from the invariance of the Haar state, (v) is a consequence of (iii) and (vi) is implied by the
formula
\[ \psi(T_{\phi}(x)) = \phi (L_{\psi}(x)),\]
valid for all $x \in \mlg$, $\psi \in \mlg_*$.
\end{proof}

All the above properties have their counterparts for the left convolution operators (this time the map $\phi \to
L_{\phi}$ is an antihomomorphism).

For $\phi \in \mlg^+_*$ we define (for each $n \in \bn$)
\begin{equation} \phi_n= \frac{1}{n} \sum_{k=1}^n {\phi^{\star k}}.
\label{avstat}\end{equation}

Properties above in conjunction with Theorem \ref{ind} imply the following fact (the notation as in the
previous subsection); the reference state on $\mlg$ will always be the Haar state.

\begin{tw} \label{inft}
 For any $\phi \in \mlg_*^+$ and $x \in \mlg$
\[ M_n (T_{\phi}) (x) = T_{\phi_n} (x) \stackrel{n \to \infty}{\longrightarrow} F(x)\]
almost uniformly.
\end{tw}

Properties of compact quantum groups allow us in fact to identify (in most of the cases) the limit in the above theorem. First
let us mention the following result due to V.\,Runde  (Corollary 3.5 in \cite{Runde}).

\begin{tw} \label{ide}
The Banach algebra $\mlg_*$ is an ideal in $\mlg^*$ (equipped with the Arens multiplication).
\end{tw}

It is elementary to check that if $\phi \in \mlg_*$, $\rho \in \mlg^*$ the Arens multiplication $\cdot$ (both left
and right version, known to coincide in this situation) may be written in terms of convolution operators:
\[ \rho \cdot \phi = \rho \circ T_{\phi}, \;\;\; \phi \cdot \rho = \rho \circ L_{\phi}.\]
Therefore the above theorem of Runde may be interpreted as the counterpart of the classical fact that for compact
groups a convolution of a bounded measure that has a density with any bounded measure gives again a measure with a
density. In two propositions below we identify the `pointwise' limits whose existence was guaranteed by theorem
\ref{inft}.

\begin{prop}
Let $\phi\in \mlg_*^+$ be a faithful state. The fixed point space of $T_{\phi}$ consists only of scalar multiples
of $1$ (in other words, $T_{\phi}$ is ergodic).
\end{prop}

\begin{proof}
Consider the restriction of $\phi$ to the w$^*$-dense compact quantum group $\alg$. As the restriction is also a
faithful state, a remark ending Section 2 of \cite{wor2} implies that for each $a \in \alg$ there is  $\phi_n(a)
\stackrel{n \to \infty}{\longrightarrow} h(a)$.  It follows (see the proof of Proposition \ref{prr}(iii)) that for
each $a \in \alg$ the sequence $(M_n (T_{\phi}) (a))_{n=1}^{\infty}$ converges to $T_h(a) = h(a) 1$  in
w$^*$-topology. Let now $\rho \in \mlg^*$ be any w$^*$-accumulation point of the sequence
$(\phi_n)_{n=1}^{\infty}$ in the unit ball of $\mlg^*$. It is easy to check that (for each $x \in \mlg$)
\[ \rho(x)= \rho (T_{\phi}(x)) = \rho (L_{\phi}(x)).\]
Theorem \ref{ide} yields normality of $\rho$, and as the first part of the proof shows that $\rho|_{\alg} =h$ and
$\alg$ is dense, we must have $\rho=h$. Therefore the projection on the fixed point space is given by the formula
$F(x) = h(x) 1$ ($x \in \mlg$).
\end{proof}

Note that in fact we did not need the theorem of Runde; it was enough to conclude by recalling the
w$^*$-continuity of $F$. The next corollary however makes essential use of Theorem \ref{ide}.

\begin{prop}
Let $\phi\in \mlg_*^+$. The sequence $(\phi_n)_{n=1}^{\infty}$ is weakly convergent to a normal functional $\rho$.
In particular, for each $x \in \mlg$
\[ M_n (T_{\phi}) (x)\stackrel{n \to \infty}{\longrightarrow} T_{\rho}(x)\]
almost uniformly.
\end{prop}

\begin{proof}
We can assume that $\phi$ is a state. Choosing this time two, potentially different, accumulation points $\rho_1,
\rho_2$ of the sequence $(\phi_n)_{n=1}^{\infty}$ in the unit ball of $\mlg^*$ we deduce as above that both
$\rho_1, \rho_2$ are normal. Theorem \ref{ind} and Properties \ref{prr} imply that in fact $T_{\rho_1} = F =
T_{\rho_2}$. Further the cancellation properties of $\alg$ yield the implication
\[T_{\rho_1}= T_{\rho_2}  \;\; \Longrightarrow \;\; \rho_1|_{\alg} = \rho_2|_{\alg},\]
and  density of $\alg$ in $\mlg$ gives the equality $\rho_1=\rho_2$.
\end{proof}

\section{Extensions to $L^p$-spaces and iterates of symmetric convolution operators}

This section will only briefly introduce bits of notation and terminology - for precise treatment of Haagerup
$L^p$-spaces we refer for example to \cite{full}.  The `density' operator of the Haar state will be denoted by $\dens$, the canonical trace-like functional on $L_1(\mlg)$ by $\tau$, $p'$ will always be the exponent conjugate to $p$. For each $\phi\in \pred$ the operator defined by
\[ T_{\phi}^{(p)} (\pdens x \pdens) = \pdens T_{\phi} (x) \pdens, \;\;\; x \in \mlg,\]
extends uniquely to a continuous operator on $L^p(\mlg)$. This follows from the fact that each $T_{\phi}$ may be
written (in a canonical way) as a linear combination of four kernels and the results of \cite{full}. One of the
main theorems of the latter paper assert the almost sure convergence of ergodic averages in $L^p$-spaces. Recall
first the definition, due to R.Jajte.

\begin{deft}
Let $ p\in [1, \infty)$, $x_n, x \in L^p (\mlg), n \in \bn$. The sequence $(x_n)_{n=1}^{\infty}$ is said to
converge almost surely (a.s.) to $x$ if for each $\epsilon > 0$ there exists a projection $e \in \mlg$ and a
family $(a_{n,k})_{n,k=1}^{\infty}$ of operators in $\mlg$ such that
\[ \phi(e^\perp) < \epsilon, \;\;\; x_n - x = \sum_{k=1}^{\infty} a_{n,k} \prdens, \;\;\; \lim_{n \to \infty} \|  \sum_{k=1}^{\infty} a_{n,k} e\| = 0.\]
Analogously the sequence $(x_n)_{n=1}^{\infty}$ is said to converge bilaterally almost surely (b.a.s.) to $x$ if
for each $\epsilon > 0$ there exists a projection $e \in \mlg$ and a family $(a_{n,k})_{n,k=1}^{\infty}$ of
operators in $\mlg$ such that
\[ \phi(e^\perp) < \epsilon, \;\;\; x_n - x = \sum_{k=1}^{\infty} \pdens a_{n,k} \pdens, \;\;\; \lim_{n \to \infty} \|  \sum_{k=1}^{\infty} e a_{n,k} e\| = 0.\]
\end{deft}

Note the following fact, which can be easily deduced from the described in the introduction properties of the
modular action (see formula \eqref{copmod}):

\begin{prop}   \label{modfact}
Let $\phi \in \mlg_*$. The operator $T_{\phi}$ commutes with the modular action of the Haar state if and only if
$\phi \circ \tau_{t} = \phi$ for each $t \in \br$.
\end{prop}

The set of all normal states satisfying the equivalent conditions formulated above will be denoted by $\modpred$.
It is easy to check that it is closed under convolution multiplication of $\pred$. Moreover the set $\modpred \cap
\starpred$ is a $^*$-subsemigroup of $\starpred$. The latter
follows from the commutation relations \eqref{coptau}-\eqref{Rscal}.

Corollary 7.12 of \cite{full} yields therefore the following theorem:

\begin{tw}
Let $\phi \in \modpred$ be a state, $x \in L^p(\mlg)$. The sequence $ (M_n(T_{\phi}^{(p)}) (x))_{n=1}^{\infty}$ is
b.a.s. (and even a.s. for $p>2$) convergent to $ F^{(p)}(x)$, where $F^{(p)}:L^p(\mlg) \to L^p(\mlg)$ denotes the
projection on the fixed points of $T_{\phi}^{(p)}$. If $\phi$ is faithful, $F^{(p)}(x) = \tau (\pprdens x)
\prdens$.
\end{tw}

Classical Stein Theorem (\cite{Stein}) and its noncommutative generalisation
(\cite{full}) allow to deduce the convergence of the iterates (as opposed to averages) of $T_{\phi}$ if it
induces a symmetric operator on the $L^2$-space. The states whose associated convolution operators satisfy this
property correspond to `symmetric' measures and can be characterised by an invariance property
with respect to the antipode. This is the context of the next proposition.

\begin{prop}
Let $\omega \in \starpred \cap \modpred$. Then $\left(T_{\omega} ^{(2)} \right)^* = T_{\omega^*} ^{(2)}$.
\end{prop}
\begin{proof}
Assume that $\omega$ is as above and $a, b \in \Tomita$. Note that Proposition \ref{modfact} implies in particular that
$T_{\omega}(a)\in \Tomita$. Moreover
\[ \lla T_{\omega}^{(2)} (\fourdens a \fourdens), \fourdens b \fourdens \rra = \tau \left( \fourdens (T_{\omega} (a))^* \fourdens \fourdens b \fourdens\right) =  \tau \left( \sigma_{\frac{i}{2}} (T_{\omega} (a)^*) b \dens \right) = \]
\[= h\left( \sigma_{\frac{i}{2}} (T_{\omega} (a)^*) b \right) = h\left( \sigma_{\frac{i}{2}} \left( (\idm \wot \ol{\omega})\Com(a^*)\right) b \right) = \]
\[= h \left((\idm \wot \ol{\omega})\Com(\sigma_{\frac{i}{2}}(a^*)) b \right) = \ol{\omega} \left( (h \wot \idm ) \Com(\sigma_{\frac{i}{2}}(a^*)) (b \ot 1) \right) =\]
\[= \ol{\omega} \circ \antip \left( (h \wot \idm ) (\sigma_{\frac{i}{2}}(a^*) \ot 1) \Com(b) \right) =
\omega^* \left( (h \wot \idm ) (\sigma_{\frac{i}{2}}(a^*) \ot 1) \Com(b) \right) = \]
\[ = (h \left( \sigma_{\frac{i}{2}}(a^*)  (\idm \wot \omega^*)\Com b \right) = \tau \left( \sigma_{\frac{i}{2}}(a^*) T_{\omega^*} (b) \dens  \right) =\]
\[= \tau \left( \twodens a^* \twodens (T_{\omega^*} (b) \right) =  \lla  \fourdens a \fourdens, T_{\omega^*}^{(2)} (\fourdens b \fourdens) \rra. \]
The claim follows now from the density of $\Tomita$ in $\mlg$.
\end{proof}

Therefore the Stein Theorem in our context implies the following result:

\begin{tw}
Let $\phi \in \starpred \cap \modpred$ be a state, $\phi = \phi^*$. For  $p\in (1, \infty)$ and $x \in L^p(\mlg)$
the sequence $ ((T_{\phi}^{(p)})^{2n} (x))_{n=1}^{\infty}$ is b.a.s. (and even a.s. for $p>2$) convergent to $
F^{(p)}(x)$, where $F^{(p)}:L^p(\mlg) \to L^p(\mlg)$ denotes the projection on the fixed points of
$(T_{\phi}^{(p)})^2$. If $x \in \mlg$ then the sequence $((T_{\phi})^{2n} (x))_{n=1}^{\infty}$ converges almost uniformly.
\end{tw}

\subsection*{Continuous semigroups}

The theorems stated above, exactly as in \cite{full}, have their multi-parameter versions and counterparts for continuous
semigroups. We mention for example the following ($F$ denotes this time a projection on the space of fixed points
of the semigroup in question):

\begin{tw}
Let $(\phi_t)_{t > 0}$ be a (weakly continuous) convolution semigroup of normal states on $\mlg$. Then for each $x
\in \mlg$
\[ M_t(x)= \frac{1}{t} \int_0^t  T_{\phi_s} (x) \textrm{ ds } \stackrel{t \to \infty}{\longrightarrow} F(x)\]
almost uniformly. If $\phi_t \in \modpred$ for all $t \geq 0$ then for every $ p \in [1, \infty)$,  $x \in L^p
(\mlg)$
\[ M^{(p)}_t(x)= \frac{1}{t} \int_0^t  T^{(p)}_{\phi_s}(x) \textrm{ ds } \stackrel{t \to \infty}{\longrightarrow} F^{(p)}(x)\]
bilaterally almost surely (and a.s. if $p >2$). If additionally $\phi_t \in \starpred \cap \modpred$, $\phi_t =
\phi_t^*$ ($t \geq 0$) then for every  $ p \in (1, \infty)$,  $x \in L^p (\mlg)$
\[  T^{(p)}_{\phi_t}(x)  \stackrel{t \to \infty}{\longrightarrow} F^{(p)}(x)\]
bilaterally almost surely (and a.s. if $p >2$).
\end{tw}


\section{Questions and comments}

The first natural question to consider is the following: what are the limit properties of the sequence
$(T_{\phi}^n = T_{\phi^{\star n}})_{n=1}^{\infty}$ if no assumption is made on symmetry properties of $\phi$? In
the classical case the general answer to this problem is given by It\^o-Kawada theorem. Suppose that $G$ is a subgroup generated of the support of the measure in question. Then the limit exists if and only if
the afore-mentioned support is not contained in a nonzero coset of any closed normal subgroup of $G$ (as otherwise a `periodicity effect' arises), and is the Haar measure on $G$
(see for example \cite{Gren}). Commutative proofs suggest that the way to obtain results of such type probably
leads through
 the Fourier analysis, which is available also for compact quantum groups.
The quantum answer is however clearly more complicated, as the example of A.\,Pal (\cite{Pal}) shows the existence
of atypical idempotent states (i.e.\ idempotent states which are not Haar measures on a quantum subgroup) on a Kac-Paljutkin
quantum group. For more examples of this type and characterisation of atypical states on various types of compact
quantum groups we refer to the forthcoming paper \cite{idem}.

The second question concerns the ergodic properties of the convolution operators on locally compact (but
noncompact) quantum groups. One difference lies in the fact that one has to deal with the left and right invariant
weights (and not states), which in general will not be equal. If discrete quantum groups are considered, the
invariant weights are strictly normal (that is, arise as  sums of normal states with orthogonal supports), as
$\mlg$ is a direct sum of matrix algebras. There is however no reason to expect that the convolution operators
would respect the underlying decomposition; their behaviour is governed by the fusion rules for unitary
(co)representations. Satisfactory general results seem to be currently out of reach, and in all probability even the
consideration of concrete examples (such as say convolution operators on the quantum deformation of the Lorentz
group) should involve the extensive use of the von Neumann algebraic techniques and exploit certain compatibility
between the modular theory of the Haar weights and the behaviour of the convolution operator in question. We hope
that the introductory results of this note may provide motivation and framework for further investigations of such
type.

\subsection*{ACKNOWLEDGMENTS}
 This paper is based on the work done during the visit of the second named
 author in Besan\c{c}on in October 2005. It arose from the discussions with
 Ren\'e Schott and Quanhua Xu. The paper was completed while the first author was
 visiting the Graduate School of Information Sciences of Tohoku University as
 Marie-Curie fellow. He would like to thank to Professors Nobuaki Obata, Fumio
 Hiai, and the other members of the GSIS for their hospitality.

\end{document}